\documentclass[12pt]{amsart}

\usepackage{amsmath, amscd, graphicx, latexsym, hyperref, times, rlepsf}

\oddsidemargin .25in \theoremstyle{plain}

\newtheorem{Thm}{Theorem}[section]
\newtheorem{Lem}[Thm]{Lemma}

\newtheorem{Prop}[Thm]{Proposition}
\newtheorem{Def}[Thm]{Definition}
\newtheorem{Rem}[Thm]{Remark}
\newtheorem{Ex}[Thm]{Example}

\newcommand{\OB}{\mathfrak{{ob}}}

\numberwithin{equation}{section}

\begin{document}

\title[]{Open books and plumbings}

\author{John B. Etnyre}

\author{Burak Ozbagci}

\address{School of Mathematics \\ Georgia Institute
of Technology \\  Atlanta  \\ Georgia}
\email{etnyre@math.gatech.edu}
\address{Department of Mathematics \\ Ko\c{c} University \\
Istanbul \\ Turkey}
\email{bozbagci@ku.edu.tr}


\date{\today}

\thanks{J.B.E. was partially supported by the NSF CAREER Grant DMS-0239600
and NSF Focused Research Grant FRG-024466.
B.O. was partially supported by the Turkish Academy of Sciences and by
the NSF Focused Research Grant FRG-024466. This research was carried out while the second author
was visiting the School of Mathematics at the Georgia Institute of
Technology. Their hospitality is gratefully acknowledged.}

\begin{abstract}
We construct, somewhat non-standard, 
Legendrian surgery diagrams for some Stein fillable contact structures on some plumbing trees of circle
bundles over spheres. We then show how to put such a surgery diagram on the pages of an open book for
$S^3,$ with relatively low genus. Thus we produce open books with low genus pages supporting these Stein fillable contact
structures, and in many cases it can be shown that these open books have minimal genus pages.
\end{abstract}

\maketitle 

\section{Introduction}

A closed and oriented 3-manifold $Y$ which is described by a
plumbing tree $\Gamma$ of oriented circle bundles over $S^2$, all
having Euler numbers less than or equal to $-2$, admits many Stein
fillable contact structures. A Kirby diagram of $Y$ is given by a
collection of unknotted circles in $S^3$ corresponding to the vertices of $\Gamma$,
linked with respect to the edges of $\Gamma$
so that the smooth framing of a circle in the diagram is exactly the
Euler number $n_i \leq -2$ corresponding to the circle bundle it
represents. To find Stein fillable contact structures on $Y$ one can
simply put all the circles into Legendrian position (with respect to
the standard contact structure in $S^3$) in such a way that the
contact framing, i.e., the Thurston-Bennequin number $tb(K_i)$ of a
circle $K_i$ is given by $n_i+1$. Then by applying Legendrian
surgery on these Legendrian circles in $S^3$ we get a Stein fillable
contact structure on $Y$. Note that, the freedom to Legendrian
realize each $K_i$ with different rotation numbers (but fixing
$tb(K_i) = n_i +1$) will enable us to find different Stein fillable
contact structures on $Y$.

We will call a plumbing tree ``non-positive" if $d_i + n_i \leq 0$
for every vertex $i$, where $d_i$ denotes the degree of the $i$-th
vertex. We will refer to a vertex in a tree with $d_i +n_i>0$ as a
bad vertex. A planar open book supporting the contact structure
obtained by a Legendrian realization of a non-positive plumbing tree
was presented in \cite {Schoenenberger05}. In this article we will generalize the
methods in \cite {Schoenenberger05} to find an open book supporting  the
contact structure obtained by a Legendrian realization of a plumbing
tree which is not necessarily non-positive.  The genus of
the open book we will construct for a tree $\Gamma$ is given by
a number $g(\Gamma)$, which we define in Section~\ref{sec:surg}.
As a preliminary step in the construction of the open books we first 
derive special Legendrian surgery diagrams for $Y$ in Section~\ref{sec:surg}.
In the following section we show how to realize these Legendrian surgery diagrams
on the pages of an open book for $S^3.$ Thus after Legendrian surgery we have an
open book supporting the desired contact structure. We also discuss how
to apply these ideas to more general contact surgery diagrams. These constructions 
lead to open books decompositions supporting all tight contact structures on small 
Seifert fibered spaces with $e_0\not = -2, -1$ having page genus zero or one. In the
last section we exhibit various examples of our construction.


\section{Open books and contact structures} \label{openbook}

Suppose that for an oriented link $L$ in a closed and oriented
3--manifold $Y$ the complement $Y\setminus L$ fibers over the
circle as $\pi \colon Y \setminus L \to S^1$ such that
$\pi^{-1}(\theta) = \Sigma_\theta $ is the interior of a compact
surface bounding $L$, for all $\theta \in S^1$. Then $(L, \pi)$ is
called an \emph{open book decomposition} (or just an \emph{open
book}) of $Y$. For each $\theta \in S^1$, the surface
$\Sigma_\theta$ is called a \emph{page}, while $L$ is called the
\emph{binding} of the open book. The monodromy of the fibration
$\pi$ is defined as the diffeomorphism of a fixed page which is
given by the first return map of a flow that is transverse to the
pages and meridional near the binding. The isotopy class of this
diffeomorphism is independent of the chosen flow and we will refer
to that as the \emph{monodromy} of the open book decomposition.

An open book $(L, \pi)$ on a 3--manifold $Y$ is said to be
\emph{isomorphic} to an open book $(L^\prime, \pi^\prime)$ on a
3--manifold $Y^\prime$, if there is a diffeomorphism $f: (Y,L) \to
(Y^\prime, L^\prime)$ such that $\pi^\prime \circ f = \pi$ on $Y
\setminus L$. In other words, an isomorphism of open books takes
binding to binding and pages to pages.

An open book can also be described as follows. First consider the
mapping torus $$\Sigma_\phi= [0,1]\times \Sigma/(1,x)\sim (0,
\phi(x))$$ where $\Sigma$ is a compact oriented surface with $r$
boundary components and $\phi$ is an element of
 the mapping class group $\Gamma_\Sigma$ of $\Sigma$.
 Since
$\phi$ is the identity map on $\partial \Sigma$,
the boundary $\partial \Sigma_\phi$ of the
mapping torus $\Sigma_\phi$
can be canonically identified with $r$ copies of $T^2 =
S^1 \times S^1$, where the first $S^1$ factor is identified with $[0,1] /
(0\sim 1)$ and the second one comes from a component of $\partial \Sigma$.
Now we glue
in $r$ copies of $D^2\times S^1$ to cap off $\Sigma_\phi$
so that $\partial D^2$ is
identified with $S^1 = [0,1] /
(0\sim 1)$ and
the $S^1$ factor in $D^2 \times S^1$
is identified with a boundary component of
$\partial \Sigma$. Thus we get a
closed $3$-manifold $Y= \Sigma_\phi \cup_{r} D^2 \times S^1 $ equipped with an open book
decomposition whose binding is the union of the
core circles $D^2 \times S^1$'s
that we glue to $\Sigma_\phi$
to obtain $Y$.
In conclusion, an element $\phi \in \Gamma_\Sigma$ determines a
$3$-manifold together with an ``abstract" open book decomposition on it.
Notice that by conjugating the monodromy $\phi$ of an open book on a 3-manifold
$Y$ by an element in $\Gamma_\Sigma$ we get
an isomorphic open book on a 3-manifold
$Y^\prime$ which is diffeomorphic to $Y$.

It has been known  for a long time that every closed and oriented
3--manifold admits an open book decomposition. Our
interest in finding open books on 3-manifolds arises from their
connection to contact structures, which we will describe very
briefly. We will assume throughout this paper that a contact
structure $\xi=\ker \alpha$ is coorientable (i.e., $\alpha$ is a
global 1--form) and positive (i.e., $\alpha \wedge d\alpha
>0$). 

{\Def \label{compatible} An open book decomposition $(L,\pi)$ of a
3--manifold $Y$ \emph{supports} a contact structure $\xi$ on $Y$ if $\xi$ can be represented by a contact form
$\alpha$ such that $\alpha ( L) > 0$ and $d \alpha > 0$ on every
page.}

\vspace{1ex}
In \cite {ThurstonWinkelnkemper75}, Thurston and Winkelnkemper show that every open book
supports a contact structure.

Suppose that an open book decomposition with page $\Sigma$ is
specified by $\phi \in \Gamma_\Sigma$. Attach a $1$-handle to the
surface $\Sigma$ connecting two points on $\partial \Sigma$ to
obtain a new surface $\Sigma^{\prime}$. Let $\gamma$ be a closed
curve in $\Sigma^{\prime}$ going over the new $1$-handle exactly
once. Define a new open book decomposition with $ \phi^\prime= \phi
\circ t_\gamma \in \Gamma_{\Sigma^{\prime}} $, where $t_\gamma$
denotes the right-handed Dehn twist along $\gamma$. The resulting
open book decomposition is called a \emph{positive stabilization} of
the one defined by $\phi$. If we use a left-handed Dehn twist
instead then we call the result a \emph{negative stabilization}. The
inverse of the above process is \index{stabilization!negative}
called \emph{positive} (\emph{negative}) \emph{destabilization}.
Notice that although the resulting monodromy depends on the chosen
curve $\gamma$, the 3--manifold specified by $(\Sigma^\prime,
\phi^\prime)$ is diffeomorphic to the 3--manifold specified by
$(\Sigma, \phi)$.

\vspace{1ex}

A converse to the Thurston-Winkelnkemper result is given by

{\Thm [Giroux \cite {Giroux02}] \label{giroux} Every contact 3--manifold
is supported by an open book. Two open books supporting 
the same contact structure admit a common positive stabilization.
Moreover two contact structures
supported by the same open book are isotopic.}

\vspace{1ex}

We refer the reader to \cite {EtnyreOBN} and \cite {OzbagciStipsicz04} for more on the
correspondence between open books and contact structures.

\section{Legendrian surgeries and plumbings}\label{sec:surg}

We assume that all the circle bundles we consider are oriented with
Euler numbers less than or equal to $-2$.  We will call a plumbing
tree of circle bundles over $S^2$ \emph{non-positive} if the sum of
the degree of the vertex and the Euler number of the bundle
corresponding to that vertex is non-positive for every vertex of the
tree. In this section we describe Legendrian surgery diagrams of 
some contact structures on plumbings of circle bundles over $S^2$ according
to trees which are not necessarily non-positive. These surgery diagrams will
be transformed into open books in the following section. 

Let us denote a circle bundle over $S^2$ with Euler number $n$ by
$Y_n$. Given a plumbing tree $\Gamma$ of circle bundles $Y_{n_i},$
denote the boundary of the plumbed sphere bundles by $Y_\Gamma.$  A
vertex with $ n_i + d_i
> 0 $ will be called a bad vertex, where
$d_i$ denote the degree (or the valence) of that vertex. We will
call a connected linear subtree $\widehat{\Gamma}\subset\Gamma$
maximal if there is no connected linear subtree
$\widetilde{\Gamma}\subset\Gamma$ such that $\widehat{\Gamma}$ is a
proper subset of $ \widetilde{\Gamma}$. The set $\Gamma \setminus
\widehat{\Gamma}$ will denote the subtree where we remove from
$\Gamma$ all the edges emanating from any vertex in
$\widehat{\Gamma} $ as well as all the vertices of
$\widehat{\Gamma}$. Take a maximal linear subtree
$\Gamma_1\subset\Gamma$ which includes at least one bad vertex. Then
take a maximal linear subtree $\Gamma_2\subset\Gamma \setminus
\Gamma_1$ which includes at last one bad vertex of $\Gamma.$ It is clear that by
iterating this process we will end up with a subtree $\Gamma
\setminus \bigcup_{j=1}^{s} \Gamma_j\subset\Gamma$ without any bad
vertices, for some disjoint subtrees $\Gamma_1, \ldots, \Gamma_s$
such that $\Gamma_{j+1} \subset \Gamma \setminus \bigcup_{t=1}^{j}
\Gamma_t\;$, for $j=1,\ldots,s-1$. Note, however, that $\Gamma_1,
\ldots, \Gamma_s$ may not be uniquely determined by $\Gamma$. In
particular, given any tree $\Gamma$, the number $s$ above is not
uniquely determined. Nevertheless there is certainly a minimum $s$,
associated to $\Gamma$, over all possible choices of subtrees in the
above process. We will refer to this number as the \emph{genus} of
$\Gamma$ and denote it by $g(\Gamma)$. If there is no bad vertex in
$\Gamma$ then we define $g(\Gamma)$ to be zero. 

{\Prop \label{surg} Suppose that we are given a plumbing tree
$\Gamma$ of $l$ circle bundles $Y_{n_i}$ such that $n_i \leq -2 $ for
all $i$. There are $(|n_1|-1)(|n_2|-1)\cdots (|n_l|-1)$ special Legendrian surgery diagrams
giving Stein manifolds with boundary 
$Y_\Gamma.$ These all have different $c_1$'s so the associated Stein fillable contact
structures are distinct.}
\begin{Rem}{\em
We do not claim these are all possible Stein fillable contact structures on 
$Y_\Gamma,$ but in some cases 
(like when $Y_\Gamma$ is a small Seifert fibered space with $e_0< -2$ ) we
do construct all Stein fillable (and all tight) contact structures. This follows from the
classification of tight contact structures in \cite{Wu??}.}
\end{Rem}
\begin{Rem}{\em
There are other, possibly more obvious, Legendrian surgery diagrams for these contact structures on
$Y_\Gamma,$ but the diagrams we derive here are the key to our constructions of open books in the 
next section.}
\end{Rem}
\begin{proof}
From \cite {Schoenenberger05} we recall how to ``roll up'' a linear plumbing tree $\Gamma.$ 
Let $\Gamma$ be the linear plumbing tree for $Y_{n_1},\ldots, Y_{n_k}$ where each $Y_{n_i}$
is plumbed to $Y_{n_{i-1}}$ and $Y_{n_{i+1}}, i=2,\ldots, k-1.$  See the left hand side of
Figure~\ref{rolledup}.
\begin{figure}[ht]
\begin{center}
 \relabelbox \small {
  \centerline{\epsfbox{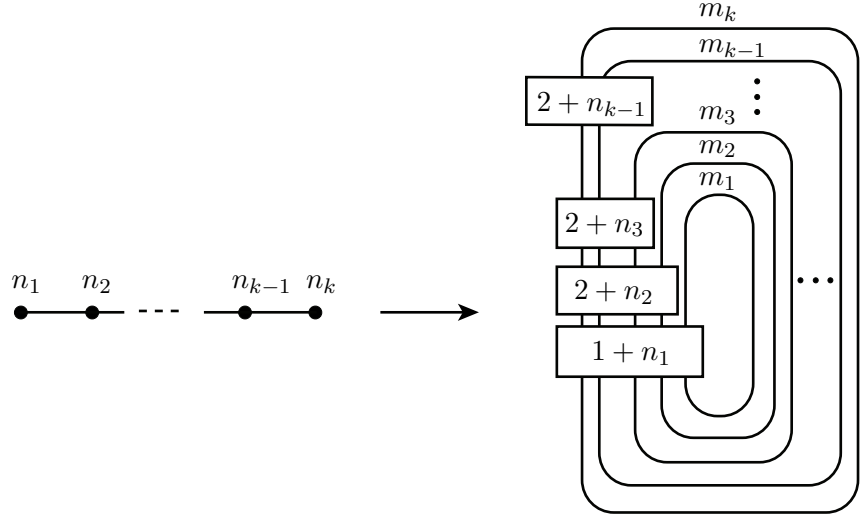}}}
  \relabel{1}{{$2+n_{k-1}$}}
  \relabel{2}{{$2+n_3$}}
  \relabel{3}{$2+n_2$}
  \relabel{4}{$1+n_1$}
  \relabel{5}{$m_k$}
  \relabel{6}{$m_{k-1}$}
  \relabel{7}{$m_{3}$}
  \relabel{8}{$m_{2}$}
  \relabel{9}{$m_1$}
  \relabel{a}{$n_1$}
  \relabel{b}{$n_2$}
  \relabel{c}{$n_{k-1}$}
  \relabel{d}{$n_k$}
  \endrelabelbox
     \caption{A linear plumbing of circle bundles and its rolled-up version.
     (The number inside a box denotes the number of full--twists we should
	apply to the knots entering into that box.)}
     \label{rolledup}
\end{center}
\end{figure}
The standard surgery diagram for $\Gamma$ is a chain of unknots $U_1, \ldots, U_k$ with
each $U_i$ simply linking $U_{i-1}$ and $U_{i+1}, i=2,\ldots, k$ and with $U_i$ having framing
$n_i.$ We think of this chain as horizontal with components labeled from left to right. Let $U_1'=U_1.$ Start with
$U_2$ and slide it over $U_1$ to get a new link with $U_2$ replaced by an unknot $U_2'$ that now
links $U_1,$ $n_1+1$ times. Now slide $U_3$ over $U_2'.$ Continue in this way until $U_k$ is slid over
$U_{k-1}'.$ The new link $L$ is called the ``rolled up" surgery diagram. See the right hand side
of Figure~\ref{rolledup}. We observe a
few salient features of this construction. First, each $U_i'$ links $U_j'$ for $j>i$ the same number of
times. Denote this linking number by $l_i.$ Secondly, $l_i\geq l_{i+1}$ for all $i.$ 
(Recall $l_i$ is negative.) In fact, $l_i=n_1+\ldots+n_{i-1} +2i-1.$ Thirdly, 
the framings $m_i$ on the $U_i'$'s are non-increasing and decrease only when $n_i<-2.$ In fact,
$m_{i+1}-m_i=n_{i+1}+2.$
Fourthly, the meridians $\mu_i$ for $U_i$ simple link $U'_i\cup \ldots \cup U'_k.$ And lastly,
$L$ sits in an unknotted solid torus neighborhood of $U_1.$
There is an obvious Legendrian representation of $L$ such that $U'_i$ is the Legendrian push off
of $U'_{i-1}$  with $|n_{i}+2|$ stabilizations. Thus Legendrian surgery produces exactly the
number of Stein manifolds claimed in the statement of the theorem.

Returning to the topological situation consider 
a tree $\Gamma$ with one valence three vertex, then we can decompose $\Gamma$
as above into linear trees $\Gamma_1$ and $\Gamma_2,$ where the first sphere bundle of 
$\Gamma_2$ is plumbed into the $i$-{th} sphere bundle of $\Gamma_1.$ Let $L_1=U_1'\cup \ldots \cup
U_k'$ and $L_2=V_1'\cup \ldots \cup V'_{k'}$ be
the rolled up surgery links for $\Gamma_1$ and $\Gamma_2,$ respectively. It is clear that if the 
neighborhood of $V_1$ in which $L_2$ sits is identified with a neighborhood of the meridian 
$\mu_i$ for $U_i$ then the resulting surgery link will describe $Y_\Gamma.$ 
As above we can Legendrian realize $L_1$ and $L_2.$ Moreover, if $n_i<-2$ then there will be
a zig-zag from the stabilization of $U_i'$ and we may link $L_2$ into $U_i'$ using this
zig-zag as shown in Figure~\ref{li}. 
\begin{figure}[ht]
  \relabelbox \small {
  \centerline{\epsfbox{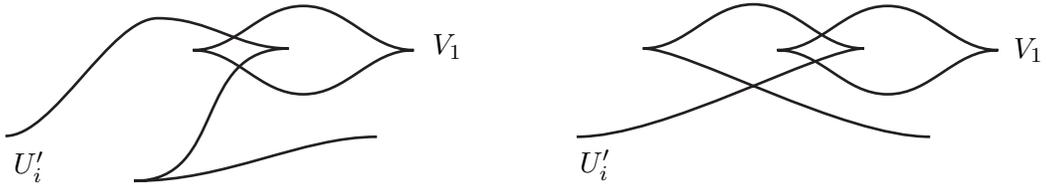}}}
  \relabel{1}{{$U_i'$}}
  \relabel{2}{{$U_i'$}}
  \relabel{3}{$V_1$}
  \relabel{4}{$V_1$}
  \endrelabelbox
        \caption{Linking $L_2$ (which is in a neighborhood fo $V_1$) to $U_i'.$}
        \label{li}
\end{figure}
If $n_i=-2$ then there is no zig-zag and no apparent way
to hook $L_2$ to $U_i'.$ However we can preform a  type 1 Legendrian Reidemeister move to 
create a cusp edge that can be used to hook $L_2$ to $U_i,$ as shown in Figure~\ref{li}. 
Thus we have Legendrian realized 
$L_1\cup L_2$ and have the desired number of Stein fillings of $Y_\Gamma.$

We can continue in
this way to obtain rolled up surgery diagrams and Legendrian surgery diagrams
for any plumbing tree.  Note that we will need to add $(n_i+d_i)$ type 1 Legendrian Reidemeister moves
to each Legendrian knot corresponding to a bad vertex.
\end{proof}

\section{Open books for plumbings which are not necessarily
non-positive.} \label{bad-vertex}

Using the notation established at the beginning of Section~\ref{sec:surg}  we are 
ready to state our main result.

{\Thm \label{tree} Suppose that we are given a plumbing tree
$\Gamma$ of circle bundles $Y_{n_i}$ such that $n_i \leq -2 $ for
all $i$. Then the Legendrian realizations of  $\;\Gamma$ from Proposition~\ref{surg}
give rise
to Stein fillable contact structures that are supported by 
open books of genus $g(\Gamma)$.}
 
This theorem was proven in \cite {Schoenenberger05} for the case with no bad vertices (compare \cite{EtguOzbagci06}). 
We generalize the ideas
there for our current proof.

\begin{Rem}{\em
In \cite{Etnyre04b} it was shown that if a contact structure is filled by a symplectic 4-manifold whose
intersection pairing does not embed in a negative definite form then the contact structure cannot be
supported by a planar open book. We observe that the intersection forms of 
some plumbings can embed in negative definite
forms but the above theorem still gives an open book with genus larger than zero. For example
if a plumbing graph has one bad vertex with Euler number $-n<0$ and valence $v<2n-l,$ where $l$ is
the number of branches from the bad vertex with length greater than 1,  then the intersection form of
this plumbing can embed into a negative definite form. It would be very interesting to see if the
genus of these open books can indeed be reduced. }
\end{Rem}

\begin{Rem}{\em 
The ideas in Theorem~\ref{tree} are much more general. Given any contact surgery diagram for
a contact structure if one can embed the individual knots in the surgery link into an open book compatible with 
the tight contact structure on $S^3$ then the ideas of  ``rolling up'', ``hooking into zig-zag's'' and
``hooking into a type 1 Legendrian Reidemeister move'' can be used to construct open books for
the resulting manifolds. While algorithms for constructing open books have been known for some time,
see for example \cite{AkbulutOzbagci01}, this algorithm seems to produce much smaller genus open books. 
We demonstrate this by constructing an open book for each tight contact structure on
small Seifert fibered spaces with $e_0\not = -1, -2.$ For notation see \cite{GhigginiLiscaStipsicz06}. Planar open books
for tight contact structures on many of these examples were also constructed in \cite{LiscaStipsicz}.}
\end{Rem}

\begin{Prop}\label{sfs}
Consider the small Seifert fibered space $M=M(e_0; r_1, r_2, r_3).$ Any tight contact structure on $M$ is
supported by an open book with planar pages if $e_0\leq -3, e_0\geq 0$ or if  $e_0=-1,\frac 12 \leq
r_1, r_2< 1$ and $0<r_3< 1.$ 
\end{Prop}

\begin{proof}[Proof of Theorem~\ref{tree}]
We recall the idea used in \cite {Schoenenberger05} to find open books
supported by a contact structure obtained by a Legendrian
realization of the linear plumbing tree described in Proposition~\ref{surg}. 
(See the proof of this proposition for notation used here.)
Consider the core
circle $\gamma$ of the open book $\OB_H$ in $S^3$ given by the
negative Hopf link $H$. The page of $\OB_H$ is an annulus and its
monodromy is a right-handed Dehn twist along $\gamma$. First
Legendrian realize $\gamma$ on a page of $\OB_H$. In \cite {Etnyre04b}, it
was shown that to stabilize a Legendrian knot on a page of an open
book, in general, one can first stabilize the open book  
and then push the knot over the 1--handle which is attached
to stabilize the open book. Apply this trick to stabilize $\OB_H$, $
(\vert m_1\vert-1 )$--times, by successively attaching 1--handles
while keeping the genus of the page to be zero. As a result, by
sliding $\gamma$ over all the attached 1--handles, we can embed the
inner-most knot in the rolled-up diagram on a page of the stabilized
open book in $S^3$ as a Legendrian knot. Then by iterating this
process from innermost to outermost knot, we can find an open book
in $S^3$ which contains all the knots in the rolled-up diagram as
Legendrian knots in distinct pages. Applying Legendrian surgery on
these knots yields a Stein fillable contact structure together with
a planar open book supporting it.

In a general tree $\Gamma$ of circle bundles $Y_{n_i} \; (n_i \leq
-2)$ without bad vertices we can take a maximal linear subtree
$\Gamma_1\subset\Gamma$ to start with and apply the algorithm above
to roll it up and construct a corresponding open book. Then take a
maximal linear subtree $\Gamma_2 \subset \Gamma \setminus \Gamma_1$
splitting off at a vertex of $\Gamma_1$. Note that there is a
stabilization used in the open book for $\Gamma_1,$ at the splitting
vertex, with its core circle so that we can apply our algorithm to
find an open book for $\Gamma_2$ starting from this annulus and
extend the previous open book. This is the translation of the left hand side of Figure~\ref{li} from
a front projection to Legendrian knots on pages of open books.
It is clear that we can continue this
process to cover all the vertices in $\Gamma$. As observed in
\cite{Schoenenberger05}, this will work as long as the tree does not have any bad
vertices, since the condition $n_i + d_i \leq  0$ guarantees that
there are as many ``free" annuli in that vertex as we need to hook
in a subtree splitting off at that vertex.  It should be clear that
we will always get a planar open book as a result.

To understand the situation with bad vertices we need to translate the right hand side of Figure~\ref{li} into
Legendrian knots on pages of open books. Specifically, we need a lemma that tells us how one can
embed a type 1 Legendrian Reidemeister move into the page of an open book.
\begin{Lem}\label{stab}
Let $(\Sigma, \phi)$ be an open book supporting a contact structure $\xi$ on $M$ and $K$ an oriented
Legendrian knot on a page of the open book.
Suppose
$R=[0,1]\times [-1,1]$ is a rectangle in the page of an open book  such that $(\partial \Sigma)\cap R=
[0,1]\times\{-1, 1\}$ and $[0,1]\times\{0\}= K\cap R$ with the orientation on $K$ agreeing with
the standard orientation on $[0,1].$ 
Stabilize the open book by adding a 1--handle to $R$ such that the 1--handle connects $[0,1]\times\{1\}$
to $[0,1]\times \{-1\}$ and the new monodromy has an extra Dehn twist along $\{\frac 12\}\times [-1,1]$
union the core of the 1--handle. Call this curve $\gamma$ and orient $\gamma$ so that the orientation
on $\gamma$ and $[-1,1]$ agree. See Figure~\ref{obli}.
The homology class  $K\pm\gamma$ can be represented by an
embedded Legendrian curve $K_\pm$ on the page. The curve $K_+$ is Legendrian isotopic to
$K$ and ``corresponds'' to a type 1 Legendrian Reidemeister move. The curve $K_-$ is isotopic to
the result of stabilizing $K$ twice, once positively and once negatively.
\end{Lem}
\begin{figure}[ht]
  \relabelbox \small {
  \centerline{\epsfbox{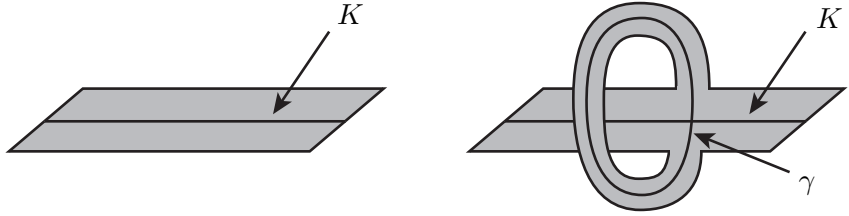}}}
  \relabel{1}{$K$}
  \relabel{2}{$K$}
  \relabel{3}{$\gamma$}
  \endrelabelbox
        \caption{The rectangle $R$ on the left and the stabilized open book on the right (the
        	embedding shown on the right is not correct there will be a full twist in the newly attached handle).}
        \label{obli}
\end{figure}
\begin{proof}
We can Legendrian realize $K\cup \gamma$ on a page of the open book. The knot $\gamma$ is a
Legendrian unknot with $tb=-1.$ Thus we can pick a disk $D$ that $\gamma$ bounds whose interior is disjoint from 
$K$ and we can make this disk convex. In the standard contact structure on $\mathbb{R}^3$ we can take an
unknot that is tangent to the $x$-axis and bounds a disk $D'$ whose interior is disjoint from the $x$-axis. 
Since $D$ and $D'$ are convex with the same dividing set (since we can assume that $D$ lies in the
complement of the binding of the original open book we know it has a tight neighborhood) we can 
assume the characteristic foliations are the same. Now we can find a contactomorphism from a 
neighborhood of $D$ union a segment of $K$ to a neighborhood of $D'$ union
a segment of the $x$-axis (so that $D$ goes to $D'$ and the segment of $K$ goes to the $x$-axis). We 
can now perform the desired operations in this local model to complete the proof.
\end{proof}

Returning to the proof of Theorem~\ref{tree},
suppose that $\Gamma$ has bad vertices. Once again we can roll up $\Gamma_1.$
When $\Gamma_1$ has bad vertices to which we wish to attach, say, $\Gamma_2$ 
we can use type 1 Legendrian Reidemeister moves as 
in the proof of Proposition~\ref{surg} to construct a Legendrian link into which $\Gamma_2$ can
be  ``hooked''.
For the open book we can stabilize as described in Lemma~\ref{stab} to create an extra annulus 
in the page of the open book that will allow us to hook in the linear graph $\Gamma_2.$
That is (using notation from the proof of Proposition~\ref{surg}) if $\Gamma_2$ is attached to $\Gamma_1$
at the unknot $U_i$ then apply Lemma~\ref{stab} to $U'_i$ (and all the subsequent $U_j'$). 
This creates an annulus in the page
of the open book that the $U_j$'s each go over exactly once for each $j\geq i.$ Now let the first unknot
of $\Gamma_2$ be a Legendrian realization of the core of the new annulus. This core will link the 
$U_j'$'s exactly once for all $j\geq i.$ We may now proceed to attach the rest of the unknots in the
rolled up version of $\Gamma_2$ as above.    
This is illustrated in Figure~\ref{badvertex}. Note that we can repeatedly apply Lemma~\ref{stab} 
to hook in arbitrarily many branches to bad  vertices of $\Gamma_1$ and only the first stabilization increases genus. Thus we see that if all the bad vertices are contained in $\Gamma_1$ the genus of
the resulting open book is one. Repeating this argument for the other $\Gamma_i$'s containing bad
vertices we see that the genus of the resulting open book is precisely $g(\Gamma).$
We have now constructed an open book for the tight contact structure on $S^3$ with the link $L$
from Proposition~\ref{surg} on its pages. Thus Legendrian surgery on this link will yield an open book
supporting the contact structure obtained by Legendrian surgery. 
\begin{figure}[ht]
\begin{center}
   \includegraphics{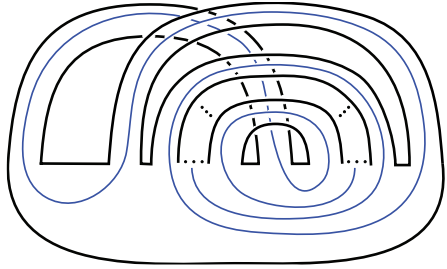}
   \caption{Positive stabilizations at a bad vertex where  $n_i + d_i > 0$.
   The number of 1--handles in the figure is given by  $1+ \vert n_i+d_i \vert$.}
	\label{badvertex}
\end{center}
\end{figure}
\end{proof}

\begin{proof}[Proof of Proposition~\ref{sfs}]
The case with $e_0\leq -3$ follows immediately from the Theorem~\ref{tree} and the classification of tight
contact structures on these manifolds from \cite{Wu??}. This result was originally proven in \cite{Schoenenberger05}. 
The case with $e_0= 0$ follows from the classification given in  \cite{GhigginiLiscaStipsicz06}. 
In particular, all these contact structures can be obtained from the tight contact structure on $S^1\times S^2$ by 
contact surgery of three Legendrian knots isotopic to $S^1\times\{pt\}.$ 
In \cite{GhigginiLiscaStipsicz06} they show that all the tight contact structures are obtained from Legendrian surgery
on Legendrian realizations of Figure~\ref{E0}.
Thus we can start with an open book for the contact structure on $S^1\times S^2$ with annular 
page and trivial monodromy. The contact framing on the components with surgery coefficient 
$a_0^j$'s is zero and since the $a_0^j$'s are all less than or equal to $-2$ we will need to stabilize the page of the
open book to Legendrian realize these components with the appropriate framing. Now if we roll up the rest of
the $a_i^j$'s onto the $a_0^j$'s, as in the proof of Proposition~\ref{surg}, 
we can easily modify the proof of Theorem~\ref{tree} to construct a genus zero open
book for these contact structures.
\begin{figure}[ht]
  \relabelbox \small {
  \centerline{\epsfbox{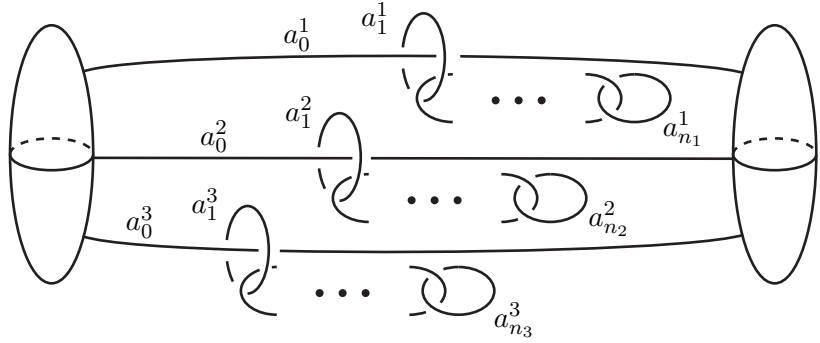}}}
  \relabel{1}{$a_0^1$}
  \relabel{2}{$a_1^1$}
  \relabel{3}{$a_{n_1}^1$}
  \relabel{4}{$a_0^2$}
  \relabel{5}{$a_1^2$}
  \relabel{6}{$a_{n_2}^2$}
  \relabel{7}{$a_0^3$}
  \relabel{8}{$a_1^3$}
  \relabel{9}{$a_{n_3}^3$}
  \endrelabelbox
        \caption{The surgery diagram for $M(0;r_1,r_2,r_3),$ here all the $a_i^j\leq-2.$}
        \label{E0}
\end{figure}
The $e_0>0$ follows similarly, the only difference with the $e_0=0$ case is that $a_0^1$ will be $-1.$ Thus
we will not be able to stabilize the knot corresponding to $a_0^1,$ 
however when one ``rolls up'' the $a_i^1$'s on $a_0^1$ we will still be able
to realize them on the page of the open book. Thus we still get a genus zero open book.
The sporadic examples with $e_0=-1$ follow from the classification of tight contact structures given
in \cite{GhigginiLiscaStipsicz??} and using the methods in the proof of Theorem~\ref{tree} to convert 
the contact structure diagrams in that paper to open books.
\end{proof}

\section{Examples}
\begin{figure}[ht]
  \relabelbox \small {
  \centerline{\epsfbox{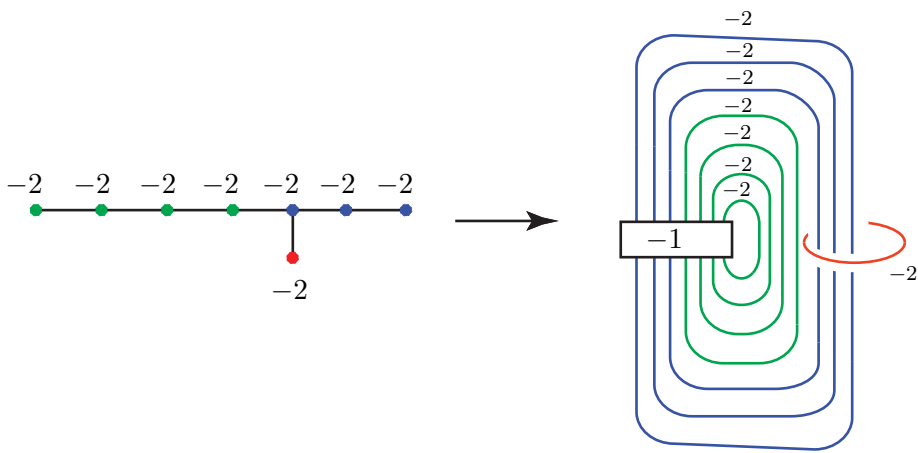}}}
  \relabel{1}{$-2$}
  \relabel{2}{$-2$}
  \relabel{3}{$-2$}
  \relabel{4}{$-2$}
  \relabel{5}{$-2$}
  \relabel{6}{$-2$}
  \relabel{7}{$-2$}
  \relabel{8}{$-2$}
  \relabel{9}{$-1$}
 \tiny \relabel{a}{$-2$}
  \relabel{b}{$-2$}
  \relabel{c}{$-2$}
  \relabel{d}{$-2$}
  \relabel{e}{$-2$}
  \relabel{f}{$-2$}
  \relabel{g}{$-2$}
  \relabel{h}{$-2$}
  \endrelabelbox
        \caption{Negative definite $E_8$-plumbing on the left and its rolled-up version on the right.}
        \label{roll235}
\end{figure}

We now demonstrate how to use the above algorithm to construct open books for various plumbing diagrams. 
{\Ex \label{one} $\,$} Consider the Poincar\'{e} homology sphere
$\Sigma(2,3,5)$ which can be given by the negative definite
$E_8$--plumbing of circle bundles over $S^2$ as in Figure~\ref{roll235}.

\begin{figure}[ht]
  \relabelbox \small {
  \centerline{\epsfbox{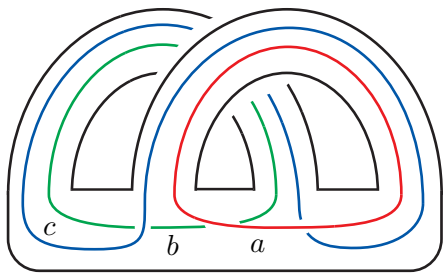}}}
  \relabel{1}{$b$}
  \relabel{2}{$a$}
  \relabel{3}{$c$}
  \endrelabelbox
  \caption{The curves $a, b$ and $c$ are embedded in \emph{distinct} pages of an open book in $S^3$ as indicated above.}
  \label{235}
\end{figure}

Also consider the genus one surface $\Sigma_{1,1}$
with one boundary component as depicted in Figure~\ref{235}. We
claim that the open book with page $\Sigma_{1,1}$ and monodromy
$$\phi = t_a^2 t_c^3 t_b^5 $$ supports the unique tight
contact structure on $\Sigma(2,3,5)$, where $t_\gamma$ denotes
a right-handed Dehn twist along a curve
$\gamma$ on a surface. Note that
$\Sigma(2,3,5)$ does not admit any planar open book supporting
its unique tight contact structure (cf. \cite {Etnyre04b}).

We now  apply our algorithm to verify the  claim about the
existence of an elliptic open book on the Poincar\'{e} homology
sphere $\Sigma(2,3,5)$, supporting  its unique tight contact
structure. The idea is to first construct an open book in $S^3$ and embed the surgery curves
into the pages of this open book so that when we perform surgeries along each of these curves 
with framing
one less than the surface framing we get $\Sigma(2,3,5)$ as the resulting
3-manifold with its associated open book.
(Note that in terms of contact structures this corresponds a Legendrian surgery.)
The monodromy of the open book in $S^3$
with page $\Sigma_{1,1}$ (as shown in Figure~\ref{235}) is $t_bt_a$.  This is obtained by
stabilizing the annulus carrying the curve $b$ in
Figure~\ref{235} by attaching the 1--handle carrying the curve $a$.
Now take the linear branch in the $E_8$-plumbing in Figure~\ref{roll235} with seven vertices.
The first four vertices (which correspond to the innermost four
curves in the rolled-up version in Figure~\ref{roll235}) are represented by parallel copies
of the curve $b$ in Figure~\ref{235}. The fifth vertex is a bad
vertex and a branch splits off with only one vertex---the eighth
vertex. The curve $c$ represents this bad vertex. The sixth and the seventh
vertices are represented by parallel copies of $c$. The eighth
vertex is represented by a parallel copy of $a$. So we embedded all the curves in
Figure~\ref{roll235} into distinct pages of the open book in $S^3$.
By performing surgeries on these curves (and taking into account
the right-handed Dehn twists we needed for the stabilizations) we
get an open book in $\Sigma(2,3,5)$ with monodromy $\phi = t_a t_c^3 t_b^4 t_b t_a$
which is equivalent to $$\phi = t_a^2 t_c^3 t_b^5 .$$ Since the monodromy is a
product of right-handed Dehn twists only, the supported contact
structure is Stein fillable (and hence tight). Therefore  this
contact structure is isotopic to the unique tight contact structure
on $\Sigma(2,3,5)$.

There is also another way of finding an elliptic open book
supporting the unique tight contact structure on
$\Sigma(2,3,5)$. The monodromy of the elliptic fibration $E(1) \to
S^2$ can be given by $(t_b t_a)^{6}$, using the notation in
Figure~\ref{235}, except that we think of the curves $a$ and $b$
embedded on a non-punctured torus. By removing the union of a
section and a cusp fiber from $E(1)$ we get a Lefschetz fibration on
the 4--manifold $W$ with punctured torus fibers whose monodromy is
$(t_bt_a)^{5}$. One can check that $\partial W$ is
diffeomorphic to $\Sigma(2,3,5)$ by Kirby calculus (see, for
example, \cite {GompfStipsicz99}). Thus there is an induced open book on
$\Sigma(2,3,5)$ with monodromy $(t_bt_a)^{5}$. Since the monodromy
of this open book is a product of right-handed Dehn twists only, the
contact structure supported by this open book is Stein fillable
(cf. \cite {Giroux02}) and in fact $W$ is a Stein filling of its boundary.
We conclude that the elliptic open book with monodromy
$(t_bt_a)^{5}$ has to support  the unique tight contact
structure on $\Sigma(2,3,5)$. Finally we note that the two elliptic
open books we described above are in fact isomorphic. In order to see the
isomorphism we first observe that $t_c = t_a^{-1} t_b t_a$. Then we plug
this relation into $t_a^2 t_c^3 t_b^5 $ to get
\begin{eqnarray*} t_a^2 t_c^3 t_b^5 &=& t_a t_b^3 t_a t_b^5 \\ &=& (t_b t_a t_b) t_b (t_b t_a t_b) t_b^3 \\ &=&
t_a t_b (t_a t_b t_a) (t_b t_a t_b) t_b^2 \\ &=& t_a t_b t_b t_a t_b t_a t_b t_a t_b^2 \\ &=&
(t_b t_a t_b) t_b t_a t_b t_a t_b t_a t_b \\ &=& t_a t_b t_a t_b t_a t_b t_a t_b t_a t_b \\ &=& (t_bt_a)^{5}.
\end{eqnarray*}
Note that we used the ``braid" relation $t_at_bt_a= t_bt_at_b$
repeatedly and cyclically permuted the words in the calculation
above.

\vspace{1ex}

\begin{figure}[ht]
  \relabelbox \small {
  \centerline{\epsfbox{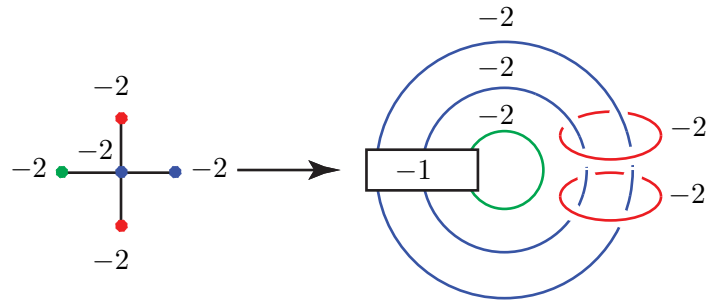}}}
  \relabel{1}{$-2$}
  \relabel{2}{$-2$}
  \relabel{3}{$-2$}
  \relabel{4}{$-2$}
  \relabel{5}{$-2$}
  \relabel{6}{$-1$}
  \relabel{7}{$-2$}
  \relabel{8}{$-2$}
  \relabel{9}{$-2$}
  \relabel{a}{$-2$}
  \relabel{b}{$-2$}
  \endrelabelbox
  \caption{A plumbing diagram on the left and its rolled-up version on the right.}
  \label{5-vertexrolledup}
\end{figure}

{\Ex \label{two} $\,$} Consider the plumbing diagram
of circle bundles and its rolled-up version shown in Figure~\ref{5-vertexrolledup}.

\begin{figure}[ht]
  \relabelbox \small {
  \centerline{\epsfbox{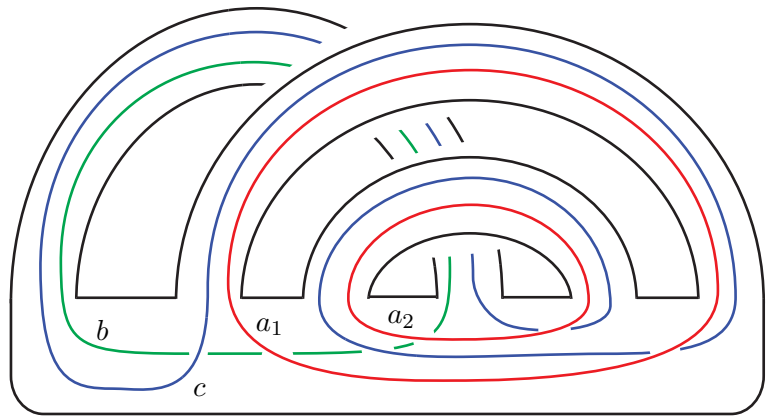}}}
  \relabel{1}{$a_1$}
  \relabel{2}{$a_2$}
  \relabel{3}{$b$}
  \relabel{4}{$c$}
  \endrelabelbox
  \caption{The curves $a_1,a_2,b$ and $c$ are embedded in \emph{distinct} pages of the open book in $S^3$ as indicated above.}
	\label{page5vertex}
\end{figure}

Then applying our algorithm we can construct an open book supporting the contact structure
obtained by the (unique) Legendrian realization of this plumbing diagram. The page $\Sigma_{1,2}$
of the open book is a torus with two
boundary components and  the monodromy is given by
$$ \phi= t_{a_1} t_{a_2} t_c^2 t_b t_b t_{a_2} t_{a_1}\;,$$ which is equivalent to the more symmetric form
$$ \phi= t_{a_1}^2 t_{a_2}^2 t_c^2 t_b^2\;,$$ where we depicted
the curves $a_1, a_2, b$ and $c$ on $\Sigma_{1,2}$ in Figure~\ref{page5vertex}. Note that Dehn twists along
two disjoint curves commute. Moreover, by plugging in $t_c= t_{a_2}^{-1} t_{a_1}^{-1} t_b t_{a_1}t_{a_2}$ we can also 
express the monodromy as  $$\phi= (t_{a_1}t_{a_2}t_b^2)^2.$$

\begin{figure}[ht]
  \relabelbox \small {
  \centerline{\epsfbox{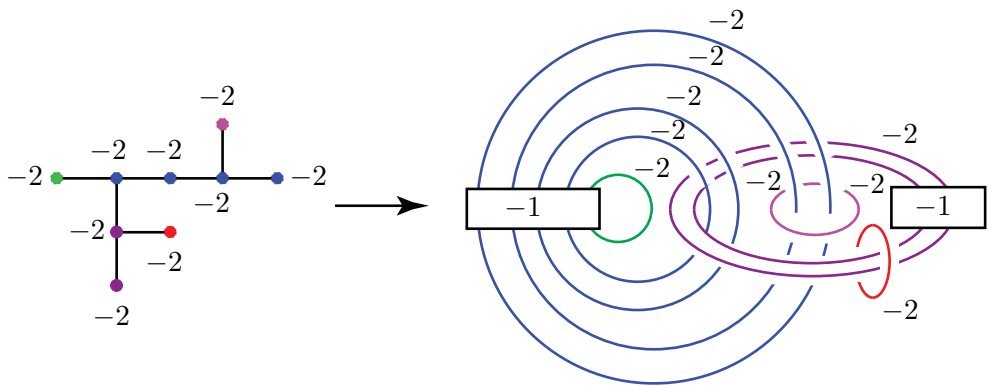}}}
  \relabel{1}{$-2$}
  \relabel{2}{$-2$}
  \relabel{3}{$-2$}
  \relabel{4}{$-2$}
  \relabel{5}{$-2$}
  \relabel{6}{$-2$}
  \relabel{7}{$-2$}
  \relabel{8}{$-2$}
  \relabel{9}{$-2$}
  \relabel{a}{$-1$}
  \relabel{b}{$-1$}
  \relabel{c}{$-2$}
  \relabel{d}{$-2$}
  \relabel{e}{$-2$}
  \relabel{f}{$-2$}
  \relabel{g}{$-2$}
  \relabel{h}{$-2$}
  \relabel{i}{$-2$}
  \relabel{j}{$-2$}
  \relabel{k}{$-2$}
  \endrelabelbox
        \caption{A plumbing diagram on the left and its rolled-up version on the right.}
        \label{ex3}
\end{figure}

{\Ex \label{three} $\,$} Consider the plumbing tree $\Gamma$ of
circle bundles and its rolled-up version shown in Figure~\ref{ex3}.

\begin{figure}[ht]
  \relabelbox \small {
  \centerline{\epsfbox{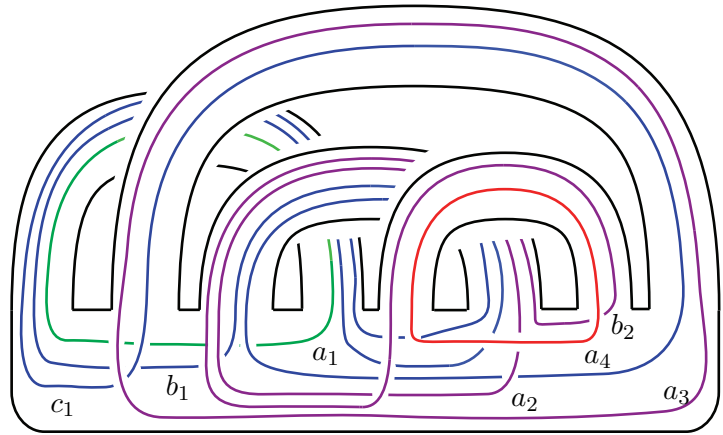}}}
  \relabel{1}{$a_1$}
  \relabel{2}{$b_1$}
  \relabel{3}{$c_1$}
  \relabel{7}{$a_2$}
  \relabel{4}{$a_4$}
  \relabel{5}{$a_3$}
  \relabel{6}{$b_2$}
  \endrelabelbox
        \caption{The curves $a_1,a_2, a_3,a_4, b_1, b_2$ and $c_1$ are embedded in \emph{distinct} pages of the open book
in $S^3$ as indicated above.}
        \label{plumb}
\end{figure}

In this final example we will illustrate how to build open books corresponding to a tree with  
two subtrees, each containing a bad vertex, meeting at a bad vertex. As the
first subtree $\Gamma_1$ take the linear tree on top with five 
vertices with a bad vertex in the second and fourth place, and as the second subtree
take $\Gamma_2$ the subtree of $\Gamma \setminus \Gamma_1$ branching from the left most bad vertex on
$\Gamma_1.$ Notice that in the
rolled-up version the part corresponding to $\Gamma_2$ is ``linked"
to the part corresponding to $\Gamma_1$. So we start with the open
book for $\Gamma_1$ and make sure the open book has been stabilized twice so that $\Gamma_2$ and 
the top most vertex can be linked into $\Gamma_1.$ Then $\Gamma_2$ is put on the pages of the open
book. To this end we must stabilize again to accommodate the bad vertex in $\Gamma_2.$
The resulting open book has page as shown in Figure~\ref{plumb}. In particular, it is a surface of genus two with one
boundary component. 
The monodromy of the open book 
supporting the contact structure obtained by Legendrian surgery on the Legendrian realization of $\Gamma$ is
then given by
$$\phi=t_{a_4}t_{a_3}t_{b_2}^2t_{c_1}^2t_{b_1}^2t_{a_1}t_{a_4}t_{a_3}t_{a_2}t_{a_1}.$$

The
careful reader might have noticed that there is in fact a maximal
linear subtree of $\Gamma$ including both bad vertices (i.e.,
$g(\Gamma)=1$) and using that as the initial subtree we could
construct a genus one open book (with two binding components) supporting the same contact
structure. However, we wanted to illustrate how to
build an open book corresponding to a tree with two subtrees, each containing a bad vertex, that meet at a bad vertex.

\def\cprime{$'$} \def\cprime{$'$}

\end{document}